%% file: NTSon_SymplecticBrockett.tex
\documentclass[hidelinks,onefignum,onetabnum]{siamart251216}

\input{ex_shared}

\ifpdf
\hypersetup{
  pdftitle={ELA_10547_R1},
  pdfauthor={N.T. Son}
}
\fi

\begin{document}

\maketitle

\begin{abstract}
The sum of symplectic eigenvalues and corresponding eigenvectors of symmetric positive-definite matrices in  
the sense of Williamson's theorem can be computed via minimization of a trace cost function under the symplecticity constraint. Optimal solutions to this problem only offer a symplectic basis for the symplectic eigenspace corresponding to the sought symplectic eigenvalues. In this note, we introduce a Brockett cost function and investigate its properties and the connection with the symplectic eigenvalues and eigenvectors of the considered matrix. Specifically, we prove that any stationary point consists of symplectic eigenvectors, characterize the saddle points and global minimizers based on which the trace minimization theorem for the symplectic eigenvalues is re-established, and the nonexistence of local nonglobal minimizers is justified.
\end{abstract}



\section{Introduction}
\label{Sec:Intro}
In his classic work \cite{Will36}, Williamson discovered that for any symmetric positive-definite (spd) matrix $M\in \Rbnn$, there exists a symplectic matrix $S\in\Rbnn$ which diagonalizes $M$ in the sense 
\begin{equation}\label{eq:sympleigdecomp}
	S^TMS = \begin{bmatrix}
		D & 0 \\ 0 & D
	\end{bmatrix},
\end{equation}
where $D=\diag(d_1,\ldots,d_n)$ with $d_j > 0, j = 1,\ldots,n.$ The decomposition \eqref{eq:sympleigdecomp} looks quite like the standard eigenvalue decomposition for real symmetric matrices 
except for the special form  on the right-hand side, which will be referred to as \emph{Williamson's diagonal form}, and the fact that $S$ is, instead of an orthogonal matrix,  a \emph{symplectic matrix}. That is, it verifies
\begin{equation}\label{eq:def_SymplMat}
	S^TJ_{2n}S = J_{2n},
\end{equation}
where $J_{2n} = \left[\begin{smallmatrix}0&I_n\\-I_n&0
\end{smallmatrix}\right]$ and $I_n$ is the $n\times n$ identity matrix. Note that the matrix $J_{2n}$ is itself symplectic, orthogonal, and skew-symmetric. It is straightforward from \eqref{eq:def_SymplMat} that $S^{-T} = J_{2n}SJ_{2n}^T$ and by multiplying this matrix to the left of \eqref{eq:sympleigdecomp}, we obtain
\begin{equation}\label{eq:sympleigdecomp_2}
	MS = J_{2n}S\begin{bmatrix}
		0 & -D \\ D & 0
	\end{bmatrix}.
\end{equation}
Let us write $S$ by columns as $S =[s_1,\ldots,s_n,s_{n+1},\ldots,s_{2n}]$. The relation \eqref{eq:sympleigdecomp_2} implies,  for any ordered vector pair $[s_j,s_{n+j}]$, that
$$M[s_j,s_{n+j}] = J_{2n}[s_j,s_{n+j}]\begin{bmatrix}
	0 & -d_j \\ d_j & 0
\end{bmatrix},\ j=1,\ldots,n,$$
which most probably suggests the name \emph{symplectic eigenvalue} for $d_j$,  \emph{symplectic eigenvector pair} associated with $d_j$ for the ordered vector pair $[s_j,s_{n+j}]$, and \emph{symplectic eigenpair} for $(d_j,[s_j,s_{n+j}])$ (cf. $Mx = I_{2n}x\lambda$ when $(\lambda,x)$ is an eigenpair of $M$.) 
More generally, denote by $\calI_k = [i_1,\ldots,i_k]$ a $k$-tuple of indices in the index set $\{1,\ldots,n\}$, and write $D_{\calI_k} = \diag(d_{i_1},\ldots,d_{i_k})$, $S_{\calI_k} = [s_{i_1},\ldots,s_{i_k},s_{n+i_1},\ldots,s_{n+i_k}]$, 
we derive from \eqref{eq:sympleigdecomp_2} the relation 
\begin{equation}\label{eq:sympleigdecomp_3}
	MS_{\calI_k} = J_{2n}S_{\calI_k}\begin{bmatrix}
		0 & -D_{\calI_k} \\ D_{\calI_k} & 0
	\end{bmatrix}.
\end{equation}
It is worth to note that the symplectic eigenvalues and eigenvectors are not the conventional ones, although some connections between them can be established 
\cite{JainM22,SonAGS21}.

Symplectic eigenvalues are involved in the study of optics \cite{EiseTRS08}, quantum mechanics \cite{KrbeTV14}, quantum computing \cite{BanchiBP2015}, quantum information/communication \cite{AdessoI2007,Weedbrooketal2012,ChakrabortyMKNC2025},  and learning linear port-Hamiltonian systems through their normal forms \cite{OrtegaY2024}. Using the connection between the symplectic eigenvalues of $M$ and the eigenvalues of the Hamiltonian matrix $J_{2n}M$, 
they can be used for stability analysis of weakly damped gyroscopic systems \cite{BennKM05}.

In principle, to compute the symplectic eigenvalues and their associated eigenvectors of $M$, one can compute the eigenvalues of the Hamiltonian matrix $J_{2n}M$ \cite{Amod06} or apply a standard eigensolver to the Hermitian pencil $M-\lambda\i J_{2n}$. However, the Hamiltonian structure, which results in special properties of its eigenvalues, 
requires careful treatment  or the computation error might lead to complex values for the Hermitian pencil \cite{BennKM05}. Moreover, some solvers do not fully deliver repeated eigenvalues. In \cite{SonAGS21}, 
inspired by the trace minimization theorem established in \cite{Hiro06,BhatJ15},
the authors addressed computing the smallest $k$ symplectic eigenvalues 
$M$, $1\leq k \leq n$, via minimizing a trace cost function under the symplecticity constraint. Namely, it seeks a solution $X_{\min}$ to 
\begin{equation*}
	\min_{X\in\R^{2n\times 2k}} \tr(X^TMX)  \quad \mbox{s.t.}\quad X^TJ_{2n}X =J_{2k}
\end{equation*}
as the first step. Here, the rectangular matrix $X$ satisfying the constraint is, similar to \eqref{eq:def_SymplMat},  also referred to as a symplectic matrix. 
In the second step, $X_{\min}^TMX_{\min}$ is diagonalized in which the obtained diagonal elements are the expected symplectic eigenvalues and the matrix product $X_{\min}X_{\diag}$, where $X_{\diag}$ diagonalizes $X_{\min}^TMX_{\min}$, is the associated symplectic eigenvectors. More steps most probably result in a larger error. One wishes to have the symplectic eigenvectors and eigenvalues right after the optimization step. Unfortunately, minimization of the trace cost function only offers a symplectic basis as the columns of $X_{\min}$ for the subspace corresponding to the smallest $k$ symplectic eigenvalues.

For eigenvalues, this issue is resolved by replacing the cost function $\tr(X^TMX)$ with the Brockett cost function $\tr(NX^TMX)$ on the set of matrices with orthogonal columns, where $N$ is a square, diagonal matrix 
\cite{Brockett1991,AbsiMS08}. More general forms were considered in 
\cite{LianWZL23,LianL2024}. 
Building on these results, in this paper, we will construct a Brockett cost function and investigate its optimization under the symplecticity constraint for the purpose of computing symplectic eigenpairs. With this in mind, we will not manage to find a general function but construct a simple one that enables our goal. After the construction of the cost function, our main results include (i) a proof that any stationary point consists of symplectic eigenvectors, (ii) an explicit formula for the global minimal value of the Brockett cost function, which allows us to re-establish the trace minimization theorem for symplectic eigenvalues mentioned above, and (iii) a justification for the nonexistence of local nonglobal minimizers of the proposed cost function.

The rest of this paper is organized as follows. After introducing notations, we  recall the necessary background on the topic in section~\ref{sec:Preliminary}. The main result will be presented in section~\ref{sec:MainResult}. Final remarks are given in section~\ref{Sec:Concl}.

\paragraph{Notation} For given integers $k, n, 1\leq k \leq n$, we denote by $\Spkn$ ($\Spn$ in case of $k=n$), $\symset(n)$, $\skewset(n)$, $\SPD(n)$ the set of $2n\times 2k$ symplectic, $n\times n$ symmetric, $n\times n$ skew-symmetric, and $n\times n$ spd matrices, respectively.   Meanwhile, $\tr(\cdot)$,  $\cdot^T$,  $\mbox{null}(\cdot)$, and $\mbox{span}(\cdot)$  
stand for the trace,  transpose, kernel of a matrix or linear operator, and subspace spanned by the column vectors (of a matrix), 
respectively. Moreover, we write $\diag(a_1,\ldots,a_n)\in\R^{n\times n}$ to mean the diagonal matrix with the entries $a_1,\ldots,a_n$ on the diagonal. 
This notation is also used for block diagonal matrices, where each $a_i$ is a submatrix block. For a~twice continuously differentiable function $f:\R^{n\times m} \to \R$, we denote by $\nabla f(X)$ and $\nabla^2 f(X)$, respectively, the Euclidean gradient and the Hessian of $f$ at~$X$. Moreover, $\mathrm{D}h(X)$ stands for the Fr\'{e}chet derivative at $X$ of a mapping $h$ between Banach spaces. Conventionally, ${\rm i}=\sqrt{-1}$ denotes the imaginary unit and $\delta_{ij}$ means the Kronecker delta.  

\section{Preliminaries}\label{sec:Preliminary}
We recall some related facts from \cite{SonAGS21}. 
The diagonalizing matrix $S$ in Williamson's decomposition \eqref{eq:sympleigdecomp} is not unique. One is referred to \cite[Thm. 3.5]{SonAGS21} for a full description of the family of diagonalizers. Thus, an argument involving $S$ or a part of it holds for any element of this family. 
A pair of column vectors $(u, v)$ in $\R^{2n}$ is said to be  \emph{symplectically normalized} if $u^T J_{2n}v = 1$. 
Two pairs of vectors $(u_1, v_1)$ and $(u_2, v_2)$  are \emph{symplectically orthogonal} if
$$u_i^TJ_{2n}v_j =  u_i^TJ_{2n}u_j = v_i^TJ_{2n}v_j  = 0 \mbox{ for } i \not = j, i, j = 1, 2.$$
A matrix $X = [u_1, \ldots, u_k, v_1, \ldots, v_k] \in  \R^{2n\times 2k}$ is symplectically normalized if each pair
$(u_i, v_i)$ is so, for $i = 1, \ldots, k$. It is  symplectically orthogonal if the pairs of
vectors $(u_i, v_i)$ are mutually symplectically orthogonal. The \emph{multiplicity} of a symplectic eigenvalue $d_j$, denoted by $m(d_j)$, is the number of times it appears in $D$. Williamson's decomposition implies that there are $m(d_j)$ symplectic eigenvector pairs, which are also symplectically normalized and orthogonal to each other, associated with $d_j$. Moreover, this set is linearly independent \cite[Thm 1.15]{deGo06} and therefore spans a $2m(d_j)$-dimensional subspace. The one associated with multiple symplectic eigenvalues is defined similarly, thanks to \cite[Cor. 2.4]{JainM22}. In view of \cite[Lem. 2.1]{SonSt22}, these subspaces are \emph{symplectic} in the sense that they, in addition to the zero vector, possess only symplectically normalized vector pairs.

A well-known result on the trace minimization for symplectic eigenvalues of an spd matrix is given in the following.
\begin{theorem}\label{theo:SymplMin}\textup{{\cite{Hiro06,BhatJ15}}}
	Let a~matrix $M\in\SPD(2n)$ have symplectic eigenvalues \mbox{$0<d_1\leq\cdots\leq d_n$}. Then for any integer $k$, $1\leq k \leq n$, it holds
	\begin{equation}\label{eq:SymplOptimProb}
		2\sum_{j=1}^{k}d_j = \min\limits_{X\in\R^{2n\times 2k}} \tr(X^TMX)   \ \mbox{ s.t. }\ X^TJ_{2n}X =J_{2k}.
	\end{equation}
\end{theorem}


Next, we collect some technical facts.
\begin{lemma}\label{lem:MainR_lem} The following statements hold.
	\begin{itemize}
		\item[i)] Given $W\in \skewset(2k)$ and $R=\diag(r_{1},\ldots,r_{2k}), r_j\not=0,$ for $j = 1,\ldots,2k$, such that
		$WR = 
		\left[\begin{smallmatrix}
			0 & -\Lambda\\
			\Lambda & \enskip 0
		\end{smallmatrix}\right]
		$,
		where $\Lambda$ is a nonsingular diagonal matrix. 
		Then $r_j = r_{j+k}$ for all $j = 1,\ldots,k$.
		\item[ii)] If a matrix $A$ commutes with a diagonal matrix $R$ with mutually distinct diagonal entries, then $A$ is also diagonal.
	\end{itemize}
\end{lemma}
\begin{proof}
	By assumption, if we write $W = \left[\begin{smallmatrix}
		W_1 & W_2\\
		-W_2^T & W_3
	\end{smallmatrix}\right]$, $W_1^T = -W_1, W_3^T = -W_3$, $R_1 = \diag(r_1,\ldots,r_k)$, and $R_2 = \diag(r_{k+1},\ldots,r_{2k})$, then we obtain immediately $W_2R_2 = R_1W_2 = -\Lambda$. Together with the nonsingularity of $R_1$ and $R_2$, these equalities first imply that $W_2$ is also a diagonal matrix. And finally, if $\Lambda$ is nonsingular, then $R_1 = R_2$.
	
	For the second statement, let us write $A = (a_{i,j})$ 
	and $R=\diag(r_1,\ldots,r_m)$ with $r_i\not=r_j$ for $i\not=j, i,j=1,\ldots,m$. The assumption in the second statement implies that $a_{i,j}(r_i-r_j)=0$ for all $i,j = 1,\ldots,m$. The conclusion is deduced from the mutual distinctness of the diagonal elements of $R$.
\end{proof}

It was pointed out in \cite[Thm. 4.6]{SonAGS21} that a solution to the minimization problem \eqref{eq:SymplOptimProb} provides a basis for the subspace associated with the smallest $k$ symplectic eigenvalues only, not a set of symplectic eigenvectors. 
The main reason is that all terms in the trace in \eqref{eq:SymplOptimProb} are equally involved. 
As a remedy to this, we consider minimizing the cost function
\begin{equation}\label{eq:Brockett}
	\min\limits_{X\in\R^{2n\times 2k}} f(X) := \tr(\tilde{N}X^TMX)
	\mbox{ s.t. } h(X):=X^TJ_{2n}X -J_{2k}=0,
\end{equation}
where $\tilde{N}$ is a diagonal matrix of positive diagonal entries. 
As \eqref{eq:Brockett} is an optimization problem under equality constraint of which the linear independence constraint qualification holds at any point in the feasible set \cite{GSAS21}, we proceed with the associated Lagrangian function
\begin{equation*}
	\calL(X,L) = \tr\left(\tilde{N}X^TMX\right) - \tr\left(L\left(X^TJ_{2n}X-J_{2k}\right)\right),
\end{equation*}
where $L\in \skewset(2k)$ is the associated Lagrange multiplier. The gradient of $\calL$ with respect to $X$ at $(X,L)$ takes the form
\begin{equation*}
	\nabla_X\calL(X,L) =  2MX\tilde{N} - 2J_{2n}XL.
\end{equation*}
Furthermore, the action of the Hessian of $\calL$ with respect to $X$ on $(Z,Z) \in \Rbnk \times \Rbnk$ reads
\begin{equation*}
	\nabla^2_{XX}\calL(X,L)[Z,Z] = 2\,\tr\left(Z^T(MZ\tilde{N} -J_{2n}Z L)\right).
\end{equation*}

Next, let us recall the first- and second-order necessary optimality conditions \cite[Prop. 3.1.1]{Bert99} for the constrained optimization problem \eqref{eq:Brockett}. A point $X_*\in \Rbnk$ is said to satisfy the \emph{first-order necessary optimality condition}, also referred to as
a~\emph{stationary point} of the problem \eqref{eq:Brockett}, if $h(X_*)=0$ and there exists 
a~Lagrange multiplier $L_* \in \skewset(2k)$ such that $\nabla_X\calL(X_*,L_*) = 0$ or equivalently $X_*^TJ_{2n}X_* = J_{2k}$ and 
\begin{equation}\label{eq:KKT}
	MX_*\tilde{N} = J_{2n}X_*L_*.
\end{equation}
The stationary point $X_*\in \Rbnk$ with the associated Lagrange multiplier $L_*$ is said to satisfy the \emph{second-order necessary optimality condition} if 
\begin{equation}\label{eq:2ndNecCond}
	\nabla^2_{XX}\calL(X_*,L_*)[Z,Z]\!=\! 2\tr\bigl(Z^T(MZ\tilde{N} -J_{2n}Z L_*)\bigr)\geq 0
\end{equation}
for all $Z \in \mbox{null}\bigl(\mathrm{D} h(X_*)\bigr)= \{Y\! \in\! \mathbb{R}^{2n\times 2k}: 
\mathrm{D} h(X_*)[Y]=Y^T\!J_{2n}X_* \!+\!X_*^TJ_{2n}Y\! = 0\}$. 

\section{Main results}\label{sec:MainResult} 
In this section, we first construct the Brockett cost function for symplectic eigenvalues and then study in detail its stationary points, minimizers, maximizers, and saddle points. 

The condition \eqref{eq:KKT} can be
equivalently written as
\begin{equation}\label{eq:KKTcond_2}
	MX_* = J_{2n}X_*L_*\tilde{N}^{-1}.
\end{equation}
We want that each stationary point is also a  symplectic eigenvector set, i.e.,   \eqref{eq:sympleigdecomp_3} holds for $X_*$ for some index subset $\calI_k$. Equating its right-hand side with that of \eqref{eq:KKTcond_2} and multiplying the derived equality from the left by $J_{2k}^TX_*^T$, 
we arrive at the condition  
$L_*\tilde{N}^{-1}=\left[\begin{smallmatrix}
	0 & -D_{\calI_k} \\ D_{\calI_k} & 0
\end{smallmatrix}\right].$ In view of Lem.~\ref{lem:MainR_lem}(i), $\tilde{N}^{-1}$ and therefore $\tilde{N}$ must take the form 
\begin{equation}\label{eq:Nform}
	\tilde{N}= \diag(N, N), N = \diag(\nu_1,\ldots,\nu_k),
\end{equation} 
with $\nu_i\not= \nu_j$ for $i\not=j,i,j=1,\ldots,k$. This observation suggests that \eqref{eq:Nform} is a necessary condition for a stationary point of \eqref{eq:Brockett} to be a symplectic eigenvector set. We will show in Prop.~\ref{theo:mainresult1} that it is also a sufficient condition. Therefore, from now on, we consider the cost function
\begin{equation}\label{eq:Brockett2}
	\min\limits_{X\in\R^{2n\times 2k}} f(X) := \tr(\tilde{N}X^TMX)
	\mbox{ s.t. } h(X):=X^TJ_{2n}X -J_{2k}=0,
\end{equation}
where $\tilde{N}$ is of the form \eqref{eq:Nform}, and mention it as the \emph{Brockett cost function} for the symplectic eigenvalues. Certainly, the order of the diagonal entries in $N$ can be arbitrary. For the sake of convenience, we will arrange them in increasing order, i.e., 
$$0< \nu_1 < \cdots < \nu_k.$$

In \cite{SonAGS21}, orthosymplectic matrices--those are both orthogonal and symplectic--play an important role in computing the symplectic eigenvalues in the optimization framework \eqref{eq:SymplOptimProb}. In this work, to treat the Brockett cost function, we introduce a specialized type of orthosymplectic matrices. Indeed, we define
\begin{equation*}
	\OrS_{\tilde{N}}(2k) = \{K\in\Rbkk : K^TJ_{2k}K = J_{2k}, K^T\tilde{N}K = \tilde{N}\}
\end{equation*}
for $\tilde{N}$ in \eqref{eq:Nform}. As a special case, it was known from \cite{SonAGS21,BabuMishra2024} that 
	\begin{equation}\label{eq:OrSp_N_form}
		\OrS_{\tilde{N}}(2k) = \bigg\{\begin{bmatrix}
			\diag(\cos(\Phi)) & \diag(\sin(\Phi)) \\
			-\diag(\sin(\Phi)) & 
			\diag(\cos(\Phi))
		\end{bmatrix}, 
		\Phi = (\phi_1,\ldots,\phi_k)\in\mathbb{R}^k\bigg\}.
	\end{equation}

\begin{remark}\label{rem:homoBrockettfunc} Using the explicit form \eqref{eq:OrSp_N_form}, it is obvious that $K\in \OrS_{\tilde{N}}(2k)$ if and only if $K^T\in \OrS_{\tilde{N}}(2k)$;   
	specifically $K\tilde{N}K^T = \tilde{N}$. Hence, it holds for any $X\in\Spkn, K\in \OrS_{\tilde{N}}(2k)$ that
	\begin{equation*}
		f(XK) = \tr(\tilde{N}K^TX^TMXK) = \tr(K\tilde{N}K^TX^TMX)
		= \tr(\tilde{N}X^TMX) = f(X).
	\end{equation*}
	That is to say, the Brockett cost function $f$ in \eqref{eq:Brockett2} is constant on the set $X\OrS_{\tilde{N}}(2k)$. 
\end{remark}

Next is a property of stationary points.
\begin{proposition}\label{lem:OrSpgroup}
	If $X_*$ is a stationary point of the minimization problem \eqref{eq:Brockett}. Then, for any $K\in\OrS_{\tilde{N}}(2k)$,  
	$X_*K$ is also a stationary point.
\end{proposition}
\begin{proof}
	First, it is direct to check that $X_*K\in\Spkn$ if $X_*\in\Spkn$ and $K\in\Spk$. Then, by assumption and Rem.~\ref{rem:homoBrockettfunc}, it holds that $$MX_*K\tilde{N}K^T=MX_*\tilde{N} = J_{2n}X_*L_* = J_{2n}X_*KK^{-1}L_*$$ 
	for some $L_*\in\skewset(2k)$. This leads to $MX_*K\tilde{N} = J_{2n}X_*KK^{-1}L_*K^{-T}$ which means that $X_*K$ is a stationary point associated with the Lagrange multiplier $K^{-1}L_*K^{-T}$.
\end{proof}

Now, we derive the expected result.
\begin{proposition}\label{theo:mainresult1}
	The matrix $X\in\Rbnk$ is a stationary point of the minimization problem \eqref{eq:Brockett2} if and only if its columns form 
	a symplectic eigenvector set of $M$.
\end{proposition}
\begin{proof}
	Assume $X_*$ is any stationary point of \eqref{eq:Brockett2} with the associated Lagrange multiplier $L_*$. From
	\eqref{eq:KKT}, we deduce that
	\begin{equation}\label{eq:MainR_prf1}
		L_* = J_{2k}^TX_*^TMX_*\tilde{N}.
	\end{equation}
	The skew-symmetry of $L_*$ immediately leads to
	\begin{equation}\label{eq:MainR_prf2}
		J_{2k}X_*^TMX_*\tilde{N} = \tilde{N}X_*^TMX_*J_{2k}^{}.
	\end{equation}
	Since $X_*^TMX_*$ is symmetric, it can  be expressed as
	\begin{equation}\label{eq:MainR_prf3}
		X_*^TMX_* = \begin{bmatrix}
			S_1&S_2\\S_2^T&S_3
		\end{bmatrix}\quad\mbox{with}\quad S_1, S_3\in \symset(k). 
	\end{equation}
	It follows from \eqref{eq:Nform},  \eqref{eq:MainR_prf2}, and \eqref{eq:MainR_prf3} that
	\begin{equation*}
		\begin{bmatrix}
			S_2^TN&S_3N\\-S_1N&-S_2N 
		\end{bmatrix} = \begin{bmatrix}
			-NS_2&NS_1\\-NS_3&NS_2^T 
		\end{bmatrix}.
	\end{equation*}
	The two identities $S_2N = -NS_2^T$ and $S_2^TN = -NS_2$ together result in $S_2N^2 = N^2S_2$. From Lem.~\ref{lem:MainR_lem}(ii), $S_2$ is diagonal which additionally yields $S_2=0$. Similarly, $S_1$ and $S_3$ are diagonal and additionally $S_1 = S_3$. Eventually, we have  
	\begin{equation}\label{eq:MainR_prf5}
		X_*^TMX_* = \begin{bmatrix}
			S_1&0\\0&S_1
		\end{bmatrix}.
	\end{equation}
	Moreover, as $M$ is spd and $X_*$ is of full rank, the diagonal elements of $S_1$ are positive. Inserting $X_*^TMX_*$ in \eqref{eq:MainR_prf5} back into \eqref{eq:MainR_prf1} and put the resulting form of $L_*$ there in \eqref{eq:KKTcond_2}, we obtain the equality 
	$$
	MX_* = J_{2n}
	X_*\begin{bmatrix}
		0&-S_1\\S_1&0
	\end{bmatrix}.
	$$
	This says $X_*$ is a symplectic eigenvector set of $M$ associated with the symplectic eigenvalues on the diagonal of $S_1$. 
	
	Conversely, assume $X_*$ is any symplectic eigenvector set associated with the symplectic 
	eigenvalues $[d_{i_1},\ldots,d_{i_k}]$. In this case \eqref{eq:sympleigdecomp_3} holds for $X_*$. Multiplying this equality with $\tilde{N}$ on the right leads to 
	$$
	MX_*\tilde{N} = J_{2n}X_*\begin{bmatrix}
		0&-D_{\calI_k}\\D_{\calI_k}&0
	\end{bmatrix}\tilde{N} = J_{2n}X_*\begin{bmatrix}
		0&-D_{\calI_k}N\\D_{\calI_k}N&0
	\end{bmatrix}.
	$$
	That is, $X_*$ is stationary point of the problem \eqref{eq:Brockett2} associated with the Lagrange multiplier $\left[\begin{smallmatrix}
		0&-D_{\calI_k}N\\D_{\calI_k}N&0
	\end{smallmatrix}\right]$. The sufficiency is then justified. 
\end{proof}

A similar result for the Brockett cost function with orthogonality constraint can be found in \cite[subsect. 4.8]{AbsiMS08} 
whose proof was mainly based on the differentiable geometry of the feasible set.

Based on Prop.~\ref{theo:mainresult1}, we will compute the global minimal value and determine the global minimizers of the cost function in \eqref{eq:Brockett2}. 
\begin{proposition}\label{theo:minimal_value}
	Assume that an spd matrix $M$ has the symplectic eigenvalues $d_1 \leq \cdots \leq d_n$, $k$ is  any number in $\{1,\ldots,n\}$, and a weighting matrix $\tilde{N}$ is of the form 
	$\diag(N, N)$, where $N = \diag(\nu_1,\ldots, \nu_k)$ with $0 < \nu_1<\nu_2<\cdots < \nu_k$. Then, it holds that
	\begin{equation}\label{eq:minBrockett}
		2(d_1\nu_k + d_2\nu_{k-1} + \cdots + d_k\nu_1)
		 = \min\limits_{X\in \Rbnk} \tr(\tilde{N}X^TMX)\ \mbox{ s.t. }\ X^TJ_{2n}X = J_{2k}.\notag
	\end{equation}
	Moreover, the columns of a minimizer $X_*$ consist of $k$ symplectic eigenvector pairs associated with the smallest $k$ symplectic eigenvalues of $M$ arranged in nonincreasing order, i.e., $X_*^TMX_*^{} = \diag(d_k,\ldots,d_1)$.
\end{proposition}
\begin{proof}
	Let us denote the columns of $X$ by $x_1,\ldots, x_{2k}$ and the diagonal elements of $\tilde{N}$ generally by $\tilde{\nu}_j, j = 1,\ldots, 2k$. 
	It holds that 
	\begin{align}\label{eq:minBrocket_proof1}
		\tr&\left(\tilde{N}X^TMX\right) = \sum_{j=1}^{2k}\tilde{\nu}_j\left(x_j^TMx_j\right)\\ \notag
		& \qquad > \nu_1\lambda_{\min}(M)\sum_{j=1}^{2k}\|x_j\|^2 = \nu_1 \lambda_{\min}(M)\|X||^2_F
	\end{align}
	for any $X\in \Rbnk$, where $\lambda_{\min}(M) >0$ is the smallest eigenvalue of $M$, $\|\cdot\|$ and $\|\cdot\|_F$ are respectively the Euclidean and Frobenius norms. This estimate, together with the unboundedness of the set $\Spkn$, implies that the Brockett cost function in \eqref{eq:Brockett2} is coercive on the closed set $\Spkn$ in $\Rbnk$. In view of Weierstra\ss's theorem, see, e.g., \cite[Prop. A.8]{Bert99}, the function in \eqref{eq:Brockett2} attains a global minimizer $X_{\min} \in \Spkn$. 
	Moreover, because of the smoothness, $X_{\min}$ must be one of the stationary points determined in Prop.~\ref{theo:mainresult1} which verifies $X_{\min}^TMX_{\min} = \left[\begin{smallmatrix}
		D_{\calI_k}&0\\0&D_{\calI_k}
	\end{smallmatrix}\right]$ for some index subset $\calI_k$. Accordingly,  a stationary value of \eqref{eq:Brockett2} can only be of the form $2\sum_{j=1}^kd_{i_j}\nu_j$ where each $d_{i_j}, j=1,\ldots,k$, is a symplectic eigenvalue of $M$ and $\nu_j$ is a diagonal element of $N$. It is the classical rearrangement inequality, see, e.g., \cite[chap. 10]{HardLP34}, that  the smallest value among all possibilities is the left-hand side of \eqref{eq:minBrockett}. 
	
	Obviously, the equality holds when $[x_j,x_{k+j}]$ is a symplectic eigenvector pair associated with the symplectic eigenvalue $d_{k-j+1}$ for $j=1,\ldots,k$. 
	The proof is complete.
\end{proof}

\begin{remark}
	Prop.~\ref{theo:minimal_value} can be obtained using \cite[Thm. 3.1]{LianL2024} with an appropriate change of variable. Their derivations are, however, different.
\end{remark}

\begin{remark}\label{rem:CorMeanValue}
	Unlike the trace cost function considered in \cite{SonAGS21}, the column order of the minimizer matters; it must be exactly associated with a $k$-tuple of symplectic eigenvalues $[d_k,\ldots,d_1]$. Therefore, even in the case $k=n$, a stationary point of the Brockett cost function in \eqref{eq:Brockett2} is not necessarily a global minimizer, cf., \cite[Cor. 4.9]{SonAGS21}. Moreover, from Rem.~\ref{rem:homoBrockettfunc}, if $X_*$ is a global minimizer, then $X_*K$ is also a global one for any $K\in\OrS_{\tilde{N}}(2k)$.
\end{remark}
\begin{remark}\label{rem:TraceMinProofAsConsequence}
	Prop.~\ref{theo:minimal_value} does not cover Thm.~\ref{theo:SymplMin} as a special case since it is required that the diagonal elements of $N$ must be mutually distinct.  However, 
	Thm.~\ref{theo:SymplMin} can be derived from Prop.~\ref{theo:minimal_value} by a simple argument. Indeed, taking $N = \diag(1+\epsilon,\ldots,1 + k\epsilon)$ for  $\epsilon >0$. Applying Prop.~\ref{theo:minimal_value} yields that
	\begin{align*}
		2&(d_1(1+k\epsilon) + d_2(1+(k-1)\epsilon) +\cdots+ d_k(1+\epsilon))\\
		&= \min\limits_{X\in \Rbnk} \tr(\tilde{N}X^TMX)\ \mbox{ s.t. }\ X^TJ_{2n}X = J_{2k}.\notag
	\end{align*}
	Thanks to the continuity of the optimal solution with respect to the parameter matrix $N$, see, e.g., \cite[Prop. 4.4]{BonnansShapiro_2000}, let $\epsilon$ in this indentity tend to zero, we obtain \eqref{eq:SymplOptimProb}.
\end{remark}

Let us turn our attention toward maximal values. First, the cost function in \eqref{eq:Brockett2} is unbounded from above due to \eqref{eq:minBrocket_proof1} and therefore it has no global maximizer. Moreover, it has no local maximizer. Verification of this statement can be done similarly to \cite[Rem. 4.13]{SonAGS21}. Indeed, a maximizer $X^*$ must be a stationary point, and in view of Prop.~\ref{theo:mainresult1}, $X^* = S_{\calI_k}$
for some index set $\calI_k$. Now, let us consider the curve $W(t) = \left[\begin{smallmatrix}
	tI_k&0\\0&I_k/t
\end{smallmatrix}\right], t \in (1-\varepsilon,1+\varepsilon), \varepsilon > 0$. The value of $f$ along $X^*W(t) \subset \Spkn$ is
$$f(X^*W(t)) = \left(t^2 + \frac{1}{t^2}\right)\sum_{j=1}^{k}\nu_jd_{i_j}.$$
Accordingly, the second derivative at $t=1$ is positive and therefore
violates the necessary optimality condition for a maximixer. Thus, if a stationary point is shown not to be a minimizer, we can conclude that it is a saddle point.  

Next, we proceed with the stationary points of the cost function in \eqref{eq:Brockett2} that are not included in Prop.~\ref{theo:minimal_value}.
\begin{proposition}\label{prop:class_criti}
	Given a stationary point $X_* = [x_{*1},\ldots,x_{*k},$ $x_{*k+1},\ldots,x_{*2k}]$ associated with the symplectic eigenvalues $[d_{i_1},\ldots,d_{i_k}]$ of matrix $M$. 
	\begin{itemize}
		\item[i)] If the $k$-tuple $[d_{i_1},\ldots,d_{i_k}]$ is not in a nonincreasing order, i.e., there are $\alpha < \beta \in \{1, \ldots, k\}$ such that $d_{i_\alpha} <d_{i_\beta}$, then $X_*$ is not a local minimizer.
		\item[ii)] Assume that
		\begin{equation}\label{eq:extracond_saddle}
			d_{k} < d_{k+1}
		\end{equation}
		and if the  $k$-tuple $[d_{i_1},\ldots,d_{i_k}]$ is in a nonincreasing order and different from the smallest $k$-tuple $[d_{k},\ldots,d_{1}]$. 
		Then $X_*$ is not a local minimizer.
	\end{itemize}
\end{proposition}
\begin{proof} For each case, using \cite[Prop. 3.3]{GSAS21}, we will construct a direction along which the second-order necessary optimality condition \eqref{eq:2ndNecCond} is violated.
	First,  
	let us denote by $\tilde{P}_{ij}\in \Rbkk, 1\leq i,j\leq k,$ the permutation matrix received from the identity matrix $I_{2k}$ by interchanging the $i$-- and $j$-- rows and $(k+i)$-- and $(k+j)$-- rows. A matrix of this type satisfies
	$$\tilde{P}^T_{ij} = \tilde{P}_{ij},\quad \tilde{P}_{ij}^2 = I_{2k},\mbox{ and } \tilde{P}_{ij}^TJ_{2k}\tilde{P}_{ij} = J_{2k}.$$
	It is straightforward to check that $Z_*:= X_*J_{2k}\tilde{P}_{\alpha\beta}\in \mbox{null}\bigl(\mathrm{D} h(X_*)\bigr)$. From the proof of Prop.~\ref{theo:mainresult1}, the Lagrange multiplier associated with $X_*$ has the explicit form $L_* = \left[\begin{smallmatrix}
		0&-D_{\calI_k}N\\D_{\calI_k}N&0
	\end{smallmatrix}\right]$, where $D_{\calI_k} = \diag(d_{i_1},\ldots,d_{i_k})$. Consequently, we have
	\begin{align*}
		\nabla^2_{XX}\calL(X_*,L_*)[Z_*,Z_*] &= 2\tr(Z_*^T(MZ_*\tilde{N}\!-\! J_{2n}Z_*L_*))\\
		&= 2\tr(\tilde{P}_{\alpha\beta}^T\diag(D_{\calI_k},D_{\calI_k})\tilde{P}_{\alpha\beta}\tilde{N} - J_{2k}\tilde{N})\\
		&= 4((d_{i_\beta}\nu_\alpha + d_{i_\alpha}\nu_\beta) - (d_{i_\alpha}\nu_\alpha + d_{i_\beta}\nu_\beta)) < 0,
	\end{align*}
	in which the last strict inequality is due to the rearrangement inequality for two sequences of mutually distinct elements. 
	
	For the second class of stationary points, we start by showing that there always exists a symplectic eigenpair $(d_c,[s_c,s_{n+c}])$, for some  $c \in \{1,\ldots,k\}$, such that 
	$d_c < d_{i_{k-(c-1)}}$ and $[s_c,s_{n+c}]$ does not belong to the symplectic subspace $\mbox{span}(X_*)$. Indeed, the assumptions imply that  $d_{i_1} > d_k.$ 
	Denote by $m(d_k)$ the multiplicity of $d_k$. Accordingly, it holds that $d_k= d_{k-1} = \cdots = d_{k-(m(d_k)-1)} > d_{k-(m(d_k))}$ and the associated symplectic subspace $\mathrm{SSp}(\{d_k,\ldots, d_{k-(m(d_k)-1)}\})$, i.e., the one spanned by all symplectic eigenvector pairs associated with the referred symplectic eigenvalues, does not has a common nonzero vector with  $\mbox{span}([x_{*1},x_{*k+1}])$ \cite[Cor. 2.4]{JainM22}. If $\mathrm{SSp}(\{d_k,\ldots, d_{k-(m(d_k)-1)}\}) \not\subset \mathrm{span}(X_*\backslash \{x_{*1},x_{*k+1}\})$, there exists a symplectic eigenvector pair, which can always be associated with $d_k$ due to the multiplicity, being outside 
	$\mathrm{span}(X_*\backslash \{x_{*1},x_{*n+1}\})$. Those vectors together with $d_k$ constitute the desired symplectic eigenpair. Otherwise, if \begin{equation}\label{eq:inclusion}
		\mathrm{SSp}(\{d_k,\ldots, d_{k-(m(d_k)-1)}\})\! \subset\! \mathrm{span}(X_*\backslash \{x_{*1}, x_{*k+1}\}),
	\end{equation} 
	the tuple $[d_{i_1},\ldots,d_{i_k}]$ must contain all $m(d_k)$ copies of $d_k$,  due to its nonincreasing order, and therefore it holds that $m(d_k) \leq k-1$. If $m(d_k) = k-1$, one can check that $d_1$ with an associated symplectic eigenvector pair is the expected symplectic eigenpair. In the other case, if $m(d_k) < k-1$, we consecutively apply these arguments to the next symplectic eigenvalue $d_{k-m(d_k)}$ and note that the inclusion situation \eqref{eq:inclusion} can not happen till the end because initially the set $\{d_{i_2},\ldots,d_{i_k}\}$ has one fewer element than $\{d_k,\ldots,d_1\}$. Thus, the mentioned symplectic eigenpair always exists. 
	
	Next, consider 
	$Z_*=[s_c,s_{n+c}][\mathrm{e}_{c,2k},\mathrm{e}_{k+c,2k}]^T$ with $\mathrm{e}_{c,2k} \in \R^{2k\times 1}$ representing the $c$-th unit vector and $[s_c,s_{n+c}]$ is the symplectic eigenvector pair associated with $d_c$ mentioned above. Due to the construction, $Z_*\in \Rbnk$ and satisfies $Z_*^TJ_{2n}X_* = 0$ which implies that $Z_*\in \mbox{null}\bigl(\mathrm{D} h(X_*)\bigr)$. Moreover, one can check that
	\begin{align*}
		\nabla^2_{XX}\calL(X_*,L_*)&[Z_*,Z_*] = 2\tr(Z_*^TMZ_*\tilde{N} - Z_*^TJ_{2n}Z_*L_*)\\
		&= 4(d_c\nu_{k-c+1} - d_{i_{k-(c-1)}}\nu_{k-c+1}) < 0.
	\end{align*}
	This last inequality completes the proof.
\end{proof}

With assumption \eqref{eq:extracond_saddle}, the nonexistence of the local nonglobal minimizer is confirmed.
\begin{corollary}
	Assume that the smallest $k$-tuple of symplectic eigenvalues of $M$ is separated from the remaining part of the symplectic spectrum. Then, the Brockett cost function in \eqref{eq:Brockett2} has no local nonglobal minimizers.
\end{corollary}
\section{Conclusion}\label{Sec:Concl}
We constructed a Brockett cost function possessing favorable properties that facilitate the study of the symplectic eigenvalue problem of spd matrices in the sense of Williamson's theorem via symplecticity-constrained optimization. Namely, we showed that the columns of any critical point form a set of symplectic eigenvectors and that every local minimizer of the function is also global. Moreover, we characterized the global minimizers, the saddle points, and computed the minimum value.

\bibliographystyle{siamplain}
\bibliography{references}
\end{document}

%% file: ex_shared.tex
\usepackage{xcolor}
\usepackage{lipsum}
\usepackage{amsfonts}
\usepackage{graphicx}
\usepackage{epstopdf}
\usepackage{algorithmic}
\ifpdf
  \DeclareGraphicsExtensions{.eps,.pdf,.png,.jpg}
\else
  \DeclareGraphicsExtensions{.eps}
\fi

\def\calL{{\cal L}}

\def\tr{\mathrm{tr}}
\def\i{\mathrm{i}}
\def\Rbkk{{\mathbb{R}^{2k\times 2k}}}
\def\Rbnk{{\mathbb{R}^{2n\times 2k}}}
\def\Rbnn{{\mathbb{R}^{2n\times 2n}}}
\def\Spkn{{\mathrm{Sp}(2k,2n)}}
\def\Spn{{\mathrm{Sp}(2n)}}
\def\Spk{{\mathrm{Sp}(2k)}}
\def\SPD{\mathrm{SPD}}
\def\diag{\mathrm{diag}}
\def\OrS{\mathrm{OrSp}}
\def\calI{\mathcal{I}}

\def\R{\mathbb{R}}
\newcommand{\skewset}{\mathrm{Skew}}
\newcommand{\symset}{\mathrm{Sym}}

\newsiamremark{remark}{Remark}
\newsiamremark{hypothesis}{Hypothesis}
\crefname{hypothesis}{Hypothesis}{Hypotheses}
\newsiamthm{claim}{Claim}
\newsiamremark{fact}{Fact}
\crefname{fact}{Fact}{Facts}

\headers{Brockett cost function for symplectic eigenvalues}{N.T. Son}

\title{Brockett cost function for symplectic eigenvalues\thanks{Submitted to the editors DATE.
\funding{This research was partially funded by the TNU--University of Sciences under the project code:~CS2024-TN06-04.}}}

\author{Nguyen Thanh Son\thanks{Department of Mathematics and Informatics, Thai Nguyen University of Sciences, 24118 Thai Nguyen, Vietnam (ntson@tnus.edu.vn).}
}

